\documentclass[12pt, a4paper]{article}
\usepackage[latin1]{inputenc}
\usepackage{vmargin}
\usepackage{epstopdf}
\usepackage[dvips]{graphicx}
\usepackage{amsfonts}
\usepackage{amssymb}
\usepackage{amsmath}
\usepackage{amsthm}
\usepackage{setspace}
\usepackage{multirow}
\usepackage{hyperref}
\usepackage{cases}
\usepackage[bottom]{footmisc}
\usepackage[small]{caption}
\usepackage[usenames,dvipsnames]{pstricks}

\setmarginsrb{1.2cm}{1.8cm}{1.2cm}{1.8cm}{0cm}{0cm}{0cm}{0.5cm}

\begin{document}

\newtheorem{Proposition}{Proposition}
\newtheorem{Theorem}{Theorem}
\newtheorem{Corollary}{Corollary}
\newtheorem{Definition}{Definition}
\newtheorem{Lemma}{Lemma}

\title{Existence and uniqueness result for mean field games with congestion effect on graphs\addtocounter{footnote}{-1}\thanks{The author wishes to acknowledge the helpful conversations with Yves Achdou (Université Paris-Diderot), François Delarue (Université Nice Sophia Antipolis), Jean-Michel Lasry (Université Paris-Dauphine) and Pierre-Louis Lions (Collège de France).}}
\author{Olivier Guéant\addtocounter{footnote}{0}\thanks{UFR de Math\'ematiques, Laboratoire Jacques-Louis Lions, Universit\'e Paris-Diderot. Avenue de France, 75013 Paris, France. \texttt{olivier.gueant@ann.jussieu.fr}}}
\date{}

\maketitle
\abstract{This paper presents a general existence and uniqueness result for mean field games equations on graphs ($\mathcal{G}$-MFG). In particular, our setting allows to take into account congestion effects of almost any form. These general congestion effects are particularly relevant in graphs in which the cost to move from one node to another may for instance depend on the proportion of players in both the source node and the target node. Existence is proved using a priori estimates and a fixed point argument \emph{à la} Schauder. We propose a new criterion to ensure uniqueness in the case of Hamiltonian functions with a complex (non-local) structure. This result generalizes the discrete counterpart of uniqueness results obtained in \cite{MFG2,MFG4}.}

\section*{Introduction}

Mean field games have been introduced in 2006 by J.-M. Lasry and P.-L. Lions \cite{MFG1,MFG2,MFG3} as the limit of a large class of stochastic differential games when the number of players increases toward infinity. Independently and almost simultaneously, Huang et al. \cite{huang} proposed a similar framework for games with a large number of identical players. The idea underlying the introduction of mean field games is that, in the case of a large number of players, interactions are such that each player only considers the statistical distribution of the others to make his decisions. From a mathematical point of view, considering empirical state distributions instead of all the individual states allows to consider games that were untractable before the introduction of mean field games. In the case of a finite number $M$ of players, finding a Nash-equilibrium is indeed almost impossible, as it involves a coupled system of $M$ nonlinear Hamilton-Jacobi-Bellman partial differential equations. Mean field games equations however boil down to two coupled equations: a backward Hamilton-Jacobi-(Bellman) equation for the value function and a forward Kolmogorov transport equation for the state distribution.\\

Research on mean field games has developed\footnote{see \cite{cardaliaguet} and \cite{bensoussanlivre} for a general appraisal.} since 2006 and we can divide the literature into 3 parts.\\

The first and wider part consists in papers on mean field game theory. A lot of papers are dedicated to proving existence and/or uniqueness for a particular class of mean field games. For instance, the case of linear-quadratic games has been widely studied by Bardi in \cite{bardi,bardi2} and by Bensoussan et \emph{al.} in \cite{bensoussanpapier}. Mean field games with quadratic Hamiltonian have been studied by Guéant in \cite{gueantquad}. A specific case where closed form solution are exhibited has also been studied in \cite{gueant2009reference}. A priori estimates to prove existence have been proposed by Gomes et \emph{al.} in \cite{gomes2012}.\footnote{Gomes also studied discrete (in time) mean field games in \cite{gomes2010discrete}.} The long term behavior of solutions has also been studied by Cardaliaguet et \emph{al.} for different classes of mean field games (see \cite{cardaliaguetlongterm1,cardaliaguetlongterm2}). The convergence of games with a finite number of players towards mean field games has also been an important research topic (see \cite{gomes2011} in addition to already cited papers). Recently, Carmona and Delarue proposed a very interesting probabilistic analysis of mean field games in \cite{carmona1} and provided a comparison between controlled McKean-Vlasov dynamics and mean field games in \cite{carmona2}.\\

In parallel with this first strand of research focussed on mean field game equations and their properties, a literature emerged on numerical methods to approximate the solution of mean field game equations. Achdou et \emph{al.} proposed finite different schemes in \cite{achdou2010mean, achdou2013} to approximate the solution of the coupled system of partial differential equations involved in mean field games. They also consider the case of the planning problem where the terminal condition on the value function is replaced by a terminal condition of the state distribution (see \cite{achdou2010capuzzo}). Gradient methods have been proposed by Lachapelle et \emph{al.} in \cite{lachapelle2010computation} in the case of potential games. The case of mean field games with quadratic Hamiltonian has also been considered, for which specific monotone numerical methods have been proposed in \cite{gueant2011,gueantnhm}.\\

The last part of the literature on mean field games is dedicated to applications. The first applications were shown in \cite{gueant2008mean} and \cite{ParisPrinceton}, mainly to economics. Other economic applications have been considered since then and an important example is the paper by Lucas and Moll \cite{moll}. These two economists indeed proposed a growth theory model based on mean field games. Chan and Sircar proposed a competition model for oil producers using mean field game equations in \cite{sircar}. Financial applications are also present in the literature. For instance, Carmona et \emph{al.} proposed a mean field game model of systemic risk. Lachapelle et \emph{al.} also proposed a stylized mean field game model to understand price formation on financial markets. Finally, applications to population dynamics have also been presented in papers on numerical methods (see also \cite{lachapellewolfram} for an example with congestion).\\

In the literature, mean field games have most often been considered on continuous state spaces. As in \cite{gomes2011}, we consider in this paper the opposite case of a discrete state space. More precisely, we consider mean field games on graphs and we do not focus on a particular structure for the graph, that is for the possibility to go from one node to another. The discrete counterpart of the mean field game partial differential equations are $2N$ ordinary differential equations, where $N$ is the cardinal of the state space (the number of nodes in the graph).\\

The goal of the paper is to provide existence and uniqueness of a solution to these $\mathcal{G}$-MFG equations. As far as existence is concerned, our framework is very general and allows for almost any local or non-local congestion effects. The result is obtained using a priori estimates and a Schauder fixed point argument. Our result on uniqueness is more general than the discrete counterpart of any existing result in the case of a continuous state space. The usual monotonicity property required for games with no congestion effect (see \cite{MFG3}) is indeed generalized with the additional requirement of a structural hypothesis on the Hamiltonian function of the problem, in order to deal with congestion (see also \cite{MFG4} for another result on partial differential equations in the specific case of local congestion).\\

In the first section, we introduce the framework and the hypotheses on the payoff and cost functions. $\mathcal{G}$-MFG equations are then introduced. In Section 2, we prove the existence of a $C^1$ solution to $\mathcal{G}$-MFG equations using a priori bounds obtained through a comparison principle and a fixed point argument \emph{à la} Schauder. Section 3 is dedicated to our uniqueness result for the solutions of $\mathcal{G}$-MFG equations. Under additional smoothness assumptions for the Hamiltonian function and using the same kind of method as the one proposed in \cite{MFG3,MFG4}, we prove uniqueness under a structural assumption on the Hamiltonian function.

\section{Mean field games on graphs}

\subsection{Notations}

We consider a directed graph $\mathcal{G}$ whose nodes are indexed by integers from 1 to $N$. For each node $i \in \mathcal{N}= \lbrace 1, \ldots, N\rbrace$ we introduce $\mathcal{V}(i) \subset \mathcal{N}\setminus \lbrace i\rbrace$ the set of nodes $j$ for which a directed edge exists from $i$ to $j$. The cardinal of this set is denoted $d_i$ and called the out-degree of the node $i$. Reciprocally, we denote $\mathcal{V}^{-1}(i) \subset \mathcal{N}\setminus \lbrace i\rbrace$ the set of nodes $j$ for which a directed edge exists from $j$ to $i$.\\

We suppose that there is a continuum of anonymous and identical players of size $1$. At any time, each player is located at a given node of the graph and hence $\mathcal{G}$ is the state space of our mean field game.\\

The process modeling the position of each player is a Markov chain in continuous time. Instantaneous transition probabilities at time $t$ are described by a collection of feedback control functions $\lambda_t(i,\cdot): \mathcal{V}(i) \to \mathbb{R}_+$ (for each node $i \in \mathcal{N}$). Throughout the text, we assume that the controls are in the admissible set $\mathcal{A}$ defined by:
$$\mathcal{A} = \left\lbrace (\lambda_t(i,j))_{t \in [0,T],i \in \mathcal{N}, j\in \mathcal{V}(i)} \mathrm{\; deterministic} \left| \forall i \in \mathcal{N}, \forall j \in \mathcal{V}(i), t \mapsto \lambda_t(i,j) \in L^{\infty}(0,T)\right. \right\rbrace,$$
where $T>0$ is the final time of the game.\\

At the macroscopic level, the distribution of the players on the graph $\mathcal{G}$ is described by a function $t \mapsto m(t) = (m(1,t), \ldots, m(N,t))$, where $m(i,t)$ simply stands for the proportion of the total population of players located at node $i$ at time $t$.\\

At any time $t$, each player can pay a cost to choose the values of the transition probabilities. We assume that the instantaneous cost of a player is given by $\mathcal{L}(i,(\lambda_{i,j})_{j \in \mathcal{V}(i)},m(t))$ if he is located at node $i$ and if he sets the value of $\lambda(i,j)$ to $\lambda_{i,j}$ (for all $j \in \mathcal{V}(i)$).\\

At time $t$, each player also earns a certain instantaneous payoff from its position on the graph. We denote $f(i,m(t))$ the instantaneous payoff of a player located at node $i$ at time $t$.\\

\emph{Remark 1: We assume that the functions $\mathcal{L}(i,\cdot,\cdot)$ and $f(i, \cdot)$ do not depend on $t$. Adding a time dependence does not add technical difficulties.}\\

At time $T$, each player has a terminal payoff depending on his position on the graph. We denote $g(i,m(T))$ the terminal payoff of a player located at node $i$ at time $T$.\\

The assumptions made on the functions $\mathcal{L}(i,\cdot,\cdot)$, $f(i, \cdot)$ and $g(i, \cdot)$  are the following:\\

\begin{itemize}
  \item Continuity: $\forall i \in \mathcal{N}$, the functions $\mathcal{L}(i,\cdot,\cdot) : \mathbb{R}_+^{d_i} \times \mathbb{R}_+^N \to \mathbb{R}_+$, $ f(i, \cdot) : \mathbb{R}_+^N \to \mathbb{R}$ and $ g(i, \cdot) : \mathbb{R}_+^N \to \mathbb{R}$   are continuous.\\
  \item Convexity of the cost functions: $\forall i \in \mathcal{N}, \forall \mu \in \mathbb{R}_+^N, \mathcal{L}(i,\cdot,\mu)$ is a strictly convex function.\\
  \item Asymptotic super-linearity of the cost functions: $$\forall i \in \mathcal{N}, \forall K>0, \quad \lim_{\lambda \in \mathbb{R}_+^{d_i}, |\lambda| \to +\infty} \inf_{\mu \in [0,K]^N} \frac{\mathcal{L}(i,\lambda,\mu)}{|\lambda|} = + \infty.$$
\end{itemize}

\emph{Remark 2: as far as the variable $\mu$ is concerned, there is no real need to define the functions outside of $\mathcal{P}_N=\lbrace (\mu_1, \ldots, \mu_N) \in \mathbb{R}_+^N , \sum_{i=1}^N \mu_i = 1\rbrace$ for our existence result. However, as we need to differentiate functions with respect to the each component of the state distribution for our uniqueness result, we prefer to define functions on $\mathbb{R}_+^N$.}\\

Associated to the cost functions $\mathcal{L}(i,\cdot,\cdot)$, we define the Hamiltonian functions $\mathcal{H}(i,\cdot,\cdot)$ by:
$$
 \forall i \in \mathcal{N}, (p,\mu) \in \mathbb{R}^{d_i}\times \mathbb{R}_+^N \mapsto \mathcal{H}(i,p,\mu) = \sup_{\lambda \in \mathbb{R}_+^{d_i}} \lambda\cdot p - \mathcal{L}(i,\lambda,\mu).
$$

We recall (see for instance \cite{cannarsa}) that:\\

\begin{itemize}
  \item $\forall i\in \mathcal{N}$, $\mathcal{H}(i,\cdot, \cdot)$ is a continuous function.\\
  \item $\forall i\in \mathcal{N}, \forall \mu \in \mathbb{R}_+^N$, $\mathcal{H}(i,\cdot, \mu)$ is a $C^1$ function and $\nabla_p \mathcal{H}(i,\cdot, \cdot)$ is a continuous function with: $$ \forall p \in \mathbb{R}^{d_i}, \forall \mu \in \mathbb{R}_+^N, \nabla_p \mathcal{H}(i,p,\mu) = \mathrm{argmax}_{\lambda \in \mathbb{R}_+^{d_i}} \lambda\cdot p - \mathcal{L}(i,\lambda,\mu).$$
\end{itemize}

\subsection{Mean field game equations}

Following the notations introduced above, for a given admissible control\footnote{We call $\lambda$ a control although it is an abuse of terminology since the controls consist in the values of $\lambda$.} $\lambda \in \mathcal{A}$ and a given function $m: t \in [0,T] \mapsto (m(1,t),\ldots,m(N,t)) \in \mathcal{P}_N$ we define the intertemporal payoff function $J_m : [0,T]\times\mathcal{N}\times\mathcal{A} \to \mathbb{R}$ by:

$$J_m(i,t,\lambda) = \mathbb{E}\left[\int_{t}^T \left(-\mathcal{L}(X_s,\lambda_s(X_s,\cdot)) + f(X_s,m(s))\right) ds + g\left(X_T,m(T)\right)\right] $$

for $(X_s)_{s \in [t,T]}$ a Markov chain on $\mathcal{G}$, starting from $i$ at time $t$, with instantaneous transition probabilities given by $(\lambda_s)_{s \in [t,T]}$.\\

From this, we can adapt the definition of a (symmetric) Nash-MFG equilibrium to our context:

\begin{Definition}[Nash-MFG symmetric equilibrium]
A differentiable function $m: t \in [0,T] \mapsto (m(1,t),\ldots,m(N,t)) \in \mathcal{P}_N$ is said to be a Nash-MFG equilibrium, if there exists an admissible control $\lambda \in \mathcal{A}$ such that:
$$\forall \tilde{\lambda} \in \mathcal{A}, \forall i \in \mathcal{N} , J_m(i,0,\lambda) \ge J_m(i,0,\tilde{\lambda}),$$
and
$$\forall i \in \mathcal{N}, \frac{d\ }{dt} m(i,t) = \sum_{j \in \mathcal{V}^{-1}(i)} \lambda_t(j,i) m(j,t) - \sum_{j \in \mathcal{V}(i)} \lambda_t(i,j) m(i,t)$$
In that case, $\lambda$ is called an optimal control.\\
\end{Definition}

The first equation corresponds to the absence of profitable unilateral deviation for a player, given a trajectory of the state distribution. The second equation is a coherence equation, stating that the evolution of the state distribution is coherent with the choice of the agents.\\

Now, let us define the $\mathcal{G}$-MFG equations associated to the above mean field game with an initial distribution $m^0 \in \mathcal{P}_N$:

\begin{Definition}[The $\mathcal{G}$-MFG equations]
The $\mathcal{G}$-MFG equations consist in a system of $2N$ equations, the unknown being $ t \in [0,T] \mapsto (u(1,t), \ldots, u(N,t), m(1,t), \ldots, m(N,t))$.\\

The first half of these equations are Hamilton-Jacobi equations and consist in the following system:
$$\forall i \in \mathcal{N}, \quad \frac{d\ }{dt} u(i,t) + \mathcal{H}\left(i,(u(j,t)-u(i,t))_{j \in \mathcal{V}(i)},m(1,t),\ldots,m(N,t)\right) + f(i,m(1,t),\ldots,m(N,t))=0$$
with $u(i,T) = g(i,m(1,T),\ldots,m(N,T))$.\\

The second half of these equations are forward transport equations:
$$\forall i \in \mathcal{N}, \quad \frac{d\ }{dt} m(i,t) = \sum_{j \in \mathcal{V}^{-1}(i)} \frac{\partial \mathcal{H} (j,\cdot,m(1,t),\ldots,m(N,t))}{\partial{p_i}}\left((u(k,t)-u(j,t))_{k \in \mathcal{V}(j)}\right) m(j,t)$$$$ - \sum_{j \in \mathcal{V}(i)} \frac{\partial \mathcal{H}(i,\cdot,m(1,t),\ldots,m(N,t))}{\partial{p_j}}\left((u(k,t)-u(i,t))_{k \in \mathcal{V}(i)}\right) m(i,t)$$
with $(m(1,0), \ldots, m(N,0)) = m^0$.\\
\end{Definition}

As usual in the mean field game literature, $u$ is the value function of players. In other words, we expect to have $u(i,t) = \sup_{\lambda \in \mathcal{A}} J_m(i,t,\lambda)$ for $m$ solving the above equations.\footnote{We do not prove any verification theorem. As we obtain a smooth solution in Theorem \ref{existence}, there is in fact no technical issue.}\\

In what follows, existence and uniqueness of solutions to the $\mathcal{G}$-MFG equations are studied. We first start with existence and our proof is based on a Schauder fixed-point argument and a priori estimates to obtain compactness. Then, we present a criterion to ensure uniqueness of $C^1$ solutions.\\

\section{Existence result}

For the existence result, we first start with a lemma stating that, for a fixed $m$, the $N$ Hamilton-Jacobi equations amongst the $\mathcal{G}$-MFG equations obey a comparison principle:

\begin{Lemma}[Comparison principle]
Let $m : [0,T] \to \mathcal{P}_N$ be a continuous function.
Let $u : t \in [0,T] \mapsto (u(1,t), \ldots, u(N,t))$ be a $C^1$ function that verifies:
$$\forall i \in \mathcal{N}, \quad -\frac{d\ }{dt} u(i,t) - \mathcal{H}\left(i,(u(j,t)-u(i,t))_{j \in \mathcal{V}(i)}, m(1,t),\ldots,m(N,t) \right) - f(i,m(1,t),\ldots,m(N,t))\le0$$
with $u(i,T) \le g(i,m(1,T),\ldots,m(N,T))$.\\

Let $v : t \in [0,T] \mapsto (v(1,t), \ldots, v(N,t))$ be a $C^1$ function that verifies:
$$\forall i \in \mathcal{N}, \quad -\frac{d\ }{dt} v(i,t) - \mathcal{H}\left(i,(v(j,t)-v(i,t))_{j \in \mathcal{V}(i)}, m(1,t),\ldots,m(N,t)\right) - f(i,m(1,t),\ldots,m(N,t)) \ge0$$
with $v(i,T) \ge g(i,m(1,T),\ldots,m(N,T))$.\\

Then, $\forall i \in \mathcal{N}, \forall t \in [0,T], v(i,t) \ge u(i,t)$.\\
\end{Lemma}

\textit{Proof:}\\

Let us consider for a given $\epsilon > 0$, a point $(i^*,t^*) \in [0,T]\times\mathcal{N}$ such that $$u(i^*,t^*) - v(i^*,t^*) - \epsilon (T-t^*) = \max_{(i,t) \in [0,T]\times\mathcal{N}}  u(i,t) - v(i,t) - \epsilon (T-t)$$

If $t^*\in [0,T)$, then $\left.\frac{d\ }{dt} \left(u(i^*,t) - v(i^*,t) - \epsilon (T-t)\right)\right|_{t=t^*} \le 0$. Also, by definition of $(i^*,t^*)$, $\forall j \in \mathcal{V}(i^*)$, $u(i^*,t^*) - v(i^*,t^*) \ge u(j,t^*) - v(j,t^*)$ and hence, by definition of $\mathcal{H}(i^*,\cdot,\cdot)$:
$$\mathcal{H}\left(i^*, (v(j,t^*)-v(i^*,t^*))_{j \in \mathcal{V}(i^*)},m(1,t^*),\ldots,m(N,t^*)\right)$$$$ \ge \mathcal{H}\left(i^*, (u(j,t^*)-u(i^*,t^*))_{j \in \mathcal{V}(i^*)},m(1,t^*),\ldots,m(N,t^*)\right) $$

Combining these inequalities we get:
$$ -\frac{d\ }{dt} u(i^*,t^*) - \mathcal{H}\left(i^*,(u(j,t^*)-u(i^*,t^*))_{j \in \mathcal{V}(i^*)},m(1,t^*),\ldots,m(N,t^*)\right)$$
$$-f(i^*,m(1,t^*),\ldots,m(N,t^*))$$
$$ \ge -\frac{d\ }{dt} v(i^*,t^*) - \mathcal{H}\left(i^*,(v(j,t^*)-v(i^*,t^*))_{j \in \mathcal{V}(i^*)},m(1,t^*),\ldots,m(N,t^*)\right) $$
$$-f(i^*,m(1,t^*),\ldots,m(N,t^*))+ \epsilon$$
But this is in contradiction with the hypotheses on $u$ and $v$.\\

Hence $t^*=T$ and $\max_{(i,t) \in [0,T]\times\mathcal{N}}  u(i,t) - v(i,t) - \epsilon (T-t) \le 0$ because of the assumptions on $u(i,T)$ and $v(i,T)$.\\

This being true for any $\epsilon > 0$, we have that $\max_{(i,t) \in [0,T]\times\mathcal{N}}  u(i,t) - v(i,t) \le 0$.\qed\\

This lemma allows to provide a bound to any solution $u$ of the $N$ Hamilton-Jacobi equations and this bound is then used to obtain compactness in order to apply Schauder's fixed point theorem.

\begin{Theorem}[Existence]
\label{existence}
Under the assumptions made in section 1, there exists a $C^1$ solution $(u,m)$ of the $\mathcal{G}$-MFG equations.
\end{Theorem}

\textit{Proof:}\\

Let $m : [0,T] \to \mathcal{P}_N$ be a continuous function.\\

Let then consider the solution $u : t \in [0,T] \mapsto (u(1,t), \ldots, u(N,t))$ to the Hamilton-Jacobi equations:

$$\forall i \in \mathcal{N}, \quad \frac{d\ }{dt} u(i,t) + \mathcal{H}\left(i,(u(j,t)-u(i,t))_{j \in \mathcal{V}(i)},m(1,t), \ldots, m(N,t)\right) + f(i,m(1,t),\ldots,m(N,t)) = 0$$
with $u(i,T) = g(i,m(1,t),\ldots,m(N,t))$.\\

This function $u$ is a well defined $C^1$ function with the following bound coming from the above lemma:

$$\sup_{i \in \mathcal{N}} \|u(i,\cdot)\|_{\infty} \le \sup_{i \in \mathcal{N}} \|g(i,\cdot)\|_{\infty} + T \sup_{i \in \mathcal{N}, \mu \in \mathcal{P}_N} |\mathcal{H}(i,0,\mu)| + T \sup_{i \in \mathcal{N}, \mu \in \mathcal{P}_N} |f(i,\mu)|.$$

Using this bound and the properties on $\mathcal{H}$, we can define a function $\tilde{m} : [0,T] \to \mathcal{P}_N$ by:

$$\forall i \in \mathcal{N}, \quad \frac{d\ }{dt} \tilde{m}(i,t) = \sum_{j \in \mathcal{V}^{-1}(i)} \frac{\partial \mathcal{H}(j,\cdot,m(1,t),\ldots,m(N,t))}{\partial{p_i}}\left((u(k,t)-u(j,t))_{k \in \mathcal{V}(j)}\right) \tilde{m}(j,t)$$$$ - \sum_{j \in \mathcal{V}(i)} \frac{\partial \mathcal{H}(i,\cdot,m(1,t),\ldots,m(N,t))}{\partial{p_j}}\left((u(k,t)-u(i,t))_{k \in \mathcal{V}(i)}\right) \tilde{m}(i,t)$$
with $(\tilde{m}(0,1), \ldots, \tilde{m}(0,N)) = m^0 \in \mathcal{P}_N$.\\

$\frac{d\tilde{m}}{dt}$ is bounded, the bounds depending only on the functions $f(i,\cdot)$, $g(i,\cdot)$ and $\mathcal{H}(i,\cdot,\cdot)$, $i \in \mathcal{N}$.\\

As a consequence, if we define $\Theta : m \in C([0,T],\mathcal{P}_N) \mapsto \tilde{m} \in C([0,T],\mathcal{P}_N)$, $\Theta$ is a continuous function (from classical ODE theory) with $\Theta(C([0,T],\mathcal{P}_N))$ a relatively compact set (because of Ascoli's Theorem and the uniform Lipschitz property we just obtained).\\

Hence, because $C([0,T],\mathcal{P}_N)$ is convex, by Schauder's fixed point theorem, there exists a fixed point $m$ to $\Theta$. If we then consider $u$ associated to $m$ by the Hamilton-Jacobi equations as above, $(u,m)$ is a $C^1$ solution to the $\mathcal{G}$-MFG equations.\qed\\

\section{Uniqueness}

Coming now to uniqueness, we use similar ideas as those used in \cite{MFG2,MFG3} or \cite{MFG4} in the case of PDEs to obtain a criterion close to the one obtained in \cite{MFG4} but adapted to graphs and generalized to non-local congestion.\\

Before stating the theorem, let us introduce a few notations. $\mathcal{M}_{r,c}$ is the set of real matrices with $r$ rows and $c$ columns. In the particular case of square matrices, $\mathcal{M}_{r,r}$ is denoted $\mathcal{M}_{r}$. We also say the a square matrix $M \in \mathcal{M}_r$ is positive, and we write $M \ge 0$, if $\forall y \in \mathbb{R}^r, y'My \ge 0$ where $y'$ denotes the transpose of $y$.\\

\begin{Theorem}[Uniqueness]
Assume that $g$ is such that:
$$\forall (\nu,\mu) \in \mathcal{P}_N \times \mathcal{P}_N, \sum_{i=1}^N (g(i,\nu_1,\ldots,\nu_N) - g(i,\mu_1,\ldots,\mu_N))(\nu_i - \mu_i) \ge 0 \implies \nu=\mu.$$

\noindent Assume that $f$ is such that:
$$\forall (\nu,\mu) \in \mathcal{P}_N \times \mathcal{P}_N, \sum_{i=1}^N (f(i,\nu_1,\ldots,\nu_N) - f(i,\mu_1,\ldots,\mu_N))(\nu_i - \mu_i) \ge 0 \implies \nu=\mu$$

\noindent Assume that $\mathcal{H}$ is such that $\forall i \in \mathcal{N}, \forall j \in \mathcal{V}(i), \frac{\partial \mathcal{H}}{\partial p_j}(i,\cdot,\cdot)$ is a $C^1$ function on $\mathbb{R}^{d_i}\times\mathbb{R}_+^N$.\\

\noindent Let us define $A: (q_1,\ldots,q_N,\mu) \in \prod_{i=1}^N\mathbb{R}^{d_i}\times\mathcal{P}_N \mapsto \left(\alpha_{ij}(q_i,\mu)\right)_{i,j} \in \mathcal{M}_{N}$ defined by:
$$\alpha_{ij}(q_i,\mu) = - \frac{\partial \mathcal{H}}{\partial \mu_j}(i,q_i,\mu).$$
Let us also define, $\forall i \in \mathcal{N}$, $B^i: (q_i,\mu) \in \mathbb{R}^{d_i}\times\mathcal{P}_N \mapsto \left(\beta^i_{jk}(q_i,\mu)\right)_{j,k} \in \mathcal{M}_{N,d_i}$ defined by:
$$\beta^i_{jk}(q_i,\mu) = \frac{\mu_i}{2} \frac{\partial^2 \mathcal{H}}{\partial \mu_j \partial q_{ik}}(i,q_i,\mu).$$
Let us also define, $\forall i \in \mathcal{N}$, $C^i: (q_i,\mu) \in \mathbb{R}^{d_i}\times\mathcal{P}_N \mapsto \left(\gamma^i_{jk}(q_i,\mu)\right)_{j,k} \in \mathcal{M}_{d_i,N}$ defined by:
$$\gamma^i_{jk}(q_i,\mu) = \frac{\mu_i}{2} \frac{\partial^2 \mathcal{H}}{\partial \mu_k\partial q_{ij}}(i,q_i,\mu).$$
Let us finally define, $\forall i \in \mathcal{N}$, $D^i: (q_i,\mu) \in \mathbb{R}^{d_i}\times\mathcal{P}_N \mapsto \left(\delta^i_{jk}(q_i,\mu)\right)_{j,k} \in \mathcal{M}_{d_i}$ defined by:
$$\delta^i_{jk}(q_i,\mu) = \mu_i \frac{\partial^2 \mathcal{H}}{\partial q_{ij}\partial q_{ik}}(i,q_i,\mu).$$

Assume that $\forall (q_1,\ldots,q_N,\mu) \in \prod_{i=1}^N\mathbb{R}^{d_i}\times\mathcal{P}_N$:

$$M(q_1,\ldots,q_N,\mu)
  =
 \begin{pmatrix}
  A(q_1,\ldots,q_N,\mu)    & B^1(q_1,\mu) & \cdots   & \cdots & \cdots & B^N(q_N,\mu) \\
  C^1(q_1,\mu)  & D^1(q_1,\mu) & 0        & \cdots & \cdots & 0  \\
  \vdots    & 0        & \ddots   & \ddots &        & \vdots \\
  \vdots    & \vdots   & \ddots   & \ddots & \ddots & \vdots \\
  \vdots    & \vdots   &          & \ddots & \ddots & 0 \\
  C^N(q_N,\mu)  & 0        & \cdots   & \cdots & 0      & D^N(q_N,\mu)
 \end{pmatrix} \ge 0.
$$

Then, if $(\widehat{u},\widehat{m})$ and $(\tilde{u},\tilde{m})$ are two $C^1$ solutions of the $\mathcal{G}$-MFG equations, we have $\widehat{m}=\tilde{m}$ and $\widehat{u}=\tilde{u}$.
\end{Theorem}

\textit{Proof:}\\

The proof of this result consists in computing in two different ways the value of
$$I = \int_0^T \sum_{i=1}^N \frac{d\ }{dt} \left((\widehat{u}(i,t) - \tilde{u}(i,t))(\widehat{m}(i,t) - \tilde{m}(i,t))\right) dt .$$ \\

We first know directly that $$I = \sum_{i=1}^N (g(i,\widehat{m}(T)) - g(i,\tilde{m}(T)))(\widehat{m}(i,T) - \tilde{m}(i,T)).$$

Now, differentiating the product we get:

$$I = -\int_0^T \sum_{i=1}^N (f(i,\widehat{m}(t)) - f(i,\tilde{m}(t)))(\widehat{m}(i,t) - \tilde{m}(i,t)) dt$$
$$+ \int_0^T  \sum_{i=1}^N \left(\widehat{m}(i,t)-\tilde{m}(i,t)\right) \Bigg{\lbrack} \mathcal{H}(i,(\tilde{u}(k,t)-\tilde{u}(i,t))_{k \in \mathcal{V}(i)},\tilde{m}(t)) - \mathcal{H}(i,(\widehat{u}(k,t)-\widehat{u}(i,t))_{k \in \mathcal{V}(i)}, \widehat{m}(t))\Bigg{\rbrack}dt$$
$$+ \int_0^T  \sum_{i=1}^N \left(\widehat{u}(i,t)-\tilde{u}(i,t)\right) \Bigg{\lbrack} \sum_{j \in \mathcal{V}^{-1}(i)} \frac{\mathcal{H}}{\partial p_i}(j,(\widehat{u}(k,t)-\widehat{u}(j,t))_{k \in \mathcal{V}(j)}, \widehat{m}(t)) \widehat{m}(j,t)$$
$$ - \sum_{j \in \mathcal{V}(i)} \frac{\mathcal{H}}{\partial p_j}(i,(\widehat{u}(k,t)-\widehat{u}(i,t))_{k \in \mathcal{V}(i)}, \widehat{m}(t)) \widehat{m}(i,t) - \sum_{j \in \mathcal{V}^{-1}(i)} \frac{\mathcal{H}}{\partial p_i}(j,(\tilde{u}(k,t)-\tilde{u}(j,t))_{k \in \mathcal{V}(j)}, \tilde{m}(t)) \tilde{m}(j,t)$$
$$ + \sum_{j \in \mathcal{V}(i)} \frac{\mathcal{H}}{\partial p_j}(i,(\tilde{u}(k,t)-\tilde{u}(i,t))_{k \in \mathcal{V}(i)}, \tilde{m}(t)) \tilde{m}(i,t) \Bigg{\rbrack}dt.$$

After reordering the terms we get:
$$I = -\int_0^T \sum_{i=1}^N (f(i,\widehat{m}(t)) - f(i,\tilde{m}(t)))(\widehat{m}(i,t) - \tilde{m}(i,t)) dt$$
$$+ \int_0^T  \sum_{i=1}^N \left(\widehat{m}(i,t)-\tilde{m}(i,t)\right) \Bigg{\lbrack} \mathcal{H}(i,(\tilde{u}(k,t)-\tilde{u}(i,t))_{k \in \mathcal{V}(i)},\tilde{m}(t)) - \mathcal{H}(i,(\widehat{u}(k,t)-\widehat{u}(i,t))_{k \in \mathcal{V}(i)}, \widehat{m}(t))\Bigg{\rbrack}dt$$
$$ +\int_0^T  \sum_{i=1}^N \widehat{m}(i,t) \sum_{j \in \mathcal{V}(i)} ((\widehat{u}(j,t) - \tilde{u}(j,t)) - (\widehat{u}(i,t) - \tilde{u}(i,t)) ) \frac{\partial \mathcal{H}}{\partial{p_j}}\left(i,(\widehat{u}(k,t)-\widehat{u}(i,t))_{k \in \mathcal{V}(i)}, \widehat{m}(t)\right)dt$$
$$ -\int_0^T  \sum_{i=1}^N \tilde{m}(i,t) \sum_{j \in \mathcal{V}(i)} ((\widehat{u}(j,t) - \tilde{u}(j,t)) - (\widehat{u}(i,t) - \tilde{u}(i,t)) ) \frac{\partial \mathcal{H}}{\partial{p_j}}\left(i,(\tilde{u}(k,t)-\tilde{u}(i,t))_{k \in \mathcal{V}(i)}, \tilde{m}(t)\right)dt,$$

\emph{i.e.}:

$$I = -\int_0^T \sum_{i=1}^N (f(i,\widehat{m}(t)) - f(i,\tilde{m}(t)))(\widehat{m}(i,t) - \tilde{m}(i,t)) dt + J,$$
where
$$J = \int_0^T  \sum_{i=1}^N \left(\widehat{m}(i,t)-\tilde{m}(i,t)\right) \Bigg{\lbrack} \mathcal{H}(i,(\tilde{u}(k,t)-\tilde{u}(i,t))_{k \in \mathcal{V}(i)},\tilde{m}(t)) - \mathcal{H}(i,(\widehat{u}(k,t)-\widehat{u}(i,t))_{k \in \mathcal{V}(i)}, \widehat{m}(t))\Bigg{\rbrack}dt$$
$$ +\int_0^T  \sum_{i=1}^N \sum_{j \in \mathcal{V}(i)} ((\widehat{u}(j,t) - \tilde{u}(j,t)) - (\widehat{u}(i,t) - \tilde{u}(i,t)) ) \Bigg{\lbrack} \widehat{m}(i,t) \frac{\partial \mathcal{H}}{\partial{p_j}}\left(i,(\widehat{u}(k,t)-\widehat{u}(i,t))_{k \in \mathcal{V}(i)}, \widehat{m}(t)\right)$$$$ - \tilde{m}(i,t) \frac{\partial \mathcal{H}}{\partial{p_j}}\left(i,(\tilde{u}(k,t)-\tilde{u}(i,t))_{k \in \mathcal{V}(i)}, \tilde{m}(t)\right) \Bigg{\rbrack} dt.$$

Now, if we denote $u^\theta(t) = \tilde{u}(t) + \theta (\widehat{u}(t)-\tilde{u}(t))$ and $m^\theta(t) = \tilde{m}(t) + \theta (\widehat{m}(t)-\tilde{m}(t))$, then we have:

$$J = \int_0^T  \sum_{i=1}^N \left(\widehat{m}(i,t)-\tilde{m}(i,t)\right) \int_0^1 \Bigg{\lbrack}\sum_{j \in \mathcal{V}(i)} -\frac{\partial \mathcal{H}}{\partial p_j}(i,(u^\theta(k,t)-u^\theta(i,t))_{k \in \mathcal{V}(i)},m^\theta(t))$$$$\times ((\widehat{u}(j,t) - \tilde{u}(j,t)) - (\widehat{u}(i,t) - \tilde{u}(i,t)) ) + \sum_{j=1}^N -\frac{\partial \mathcal{H}}{\partial \mu_j}(i,(u^\theta(k,t)-u^\theta(i,t))_{k \in \mathcal{V}(i)},m^\theta(t))(\widehat{m}(j,t) - \tilde{m}(j,t))\Bigg{\rbrack}   d\theta dt$$
$$ +\int_0^T  \sum_{i=1}^N \sum_{j \in \mathcal{V}(i)} \int_0^1 \Bigg{\lbrack} \left(\widehat{m}(i,t) - \tilde{m}(i,t)\right) \frac{\partial \mathcal{H}}{\partial{p_j}}\left(i,(u^\theta(k,t)-u^\theta(i,t))_{k \in \mathcal{V}(i)}, m^\theta(t)\right)$$$$ + m^\theta(i,t) \sum_{l \in \mathcal{V}(i)} ((\widehat{u}(l,t) - \tilde{u}(l,t)) - (\widehat{u}(i,t) - \tilde{u}(i,t)) ) \frac{\partial^2 \mathcal{H}}{\partial{p_l}\partial{p_j}}\left(i,(u^\theta(k,t)-u^\theta(i,t))_{k \in \mathcal{V}(i)}, m^\theta(t)\right)$$$$ + m^\theta(i,t) \sum_{l=1}^N (\widehat{m}(l,t) - \tilde{m}(l,t)) \frac{\partial^2 \mathcal{H}}{\partial \mu_l\partial{p_j}}\left(i,(u^\theta(k,t)-u^\theta(i,t))_{k \in \mathcal{V}(i)}, m^\theta(t)\right) \Bigg{\rbrack} d\theta$$$$ \times((\widehat{u}(j,t) - \tilde{u}(j,t)) - (\widehat{u}(i,t) - \tilde{u}(i,t)) ) dt.$$

We see that the first group of terms cancels with the third one and we can write $J$ as:

$$J = \int_0^T \int_0^1 V(t) M((u^\theta(k,t)-u^\theta(1,t))_{k \in \mathcal{V}(1)}, \ldots, (u^\theta(k,t)-u^\theta(N,t))_{k \in \mathcal{V}(N)},m^\theta(t)) V(t)' d\theta dt,$$
where
$$V(t) = (\widehat{m}(t)-\tilde{m}(t), ((\widehat{u}(k,t) - \tilde{u}(k,t)) - (\widehat{u}(1,t) - \tilde{u}(1,t)) )_{k \in \mathcal{V}(1)}, \ldots,$$$$ ((\widehat{u}(k,t) - \tilde{u}(k,t)) - (\widehat{u}(N,t) - \tilde{u}(N,t)) )_{k \in \mathcal{V}(N)}).$$

Hence $J \ge 0$ and we have:
$$\int_0^T \sum_{i=1}^N (f(i,\widehat{m}(t)) - f(i,\tilde{m}(t)))(\widehat{m}(i,t) - \tilde{m}(i,t)) dt$$$$ + \sum_{i=1}^N (g(i,\widehat{m}(T)) - g(i,\tilde{m}(T)))(\widehat{m}(i,T) - \tilde{m}(i,T)) \ge 0.$$

Using the hypotheses on $f$ and $g$ we get $\widehat{m}=\tilde{m}$.\\

The comparison principle stated in Lemma 1 then brings $\widehat{u}=\tilde{u}$ and the result is proved.\qed\\

\bibliographystyle{plain}

\end{document}